\documentclass{article}
\usepackage{maa-monthly}

\newtheorem*{proposition}{Proposition}
\providecommand{\Rp}{{\mathbb{R}^n_{\geq 0}}}

\begin{document}

\title{Nonnegative Eigenvectors of Symmetric Matrices}
\markright{Notes}
\author{Hunter Swan}

\maketitle

\begin{abstract}
For matrices with all nonnegative entries, the Perron--Frobenius theorem guarantees the existence of an eigenvector with all nonnegative components.  We show that the existence of such an eigenvector is also guaranteed for a very different class of matrices, namely real symmetric matrices with exactly two eigenvalues.  We also prove a partial converse, that among real symmetric matrices with any more than two eigenvalues there exist some having no nonnegative eigenvector. 

\end{abstract}

A \textit{nonnegative vector} is one whose components are all nonnegative.  This concept has no place in pure linear algebra, as it is highly basis dependent.  However, nonnegative vectors (and their cousins, positive vectors) sometimes crop up and prove useful in applications.  For example, one of the consequences of the Perron--Frobenius theorem is that a matrix with nonnegative entries has a nonnegative (or even positive, under appropriate hypotheses) eigenvector, which fact is of great consequence for, e.g., ranking pages in search engine results~\cite{Pillai}.  

In this note, we prove that the existence of a nonnegative eigenvector is also guaranteed for a very different class of matrices, namely real symmetric matrices having only two distinct eigenvalues.  Recall that a symmetric matrix has a set of orthogonal eigenvectors that span the ambient space.  This is the only fact about symmetric matrices that we will need.

Let $M\in\mathbb{R}^{n\times n}$ be our matrix of interest.  Since we suppose $M$ has only two eigenvalues, it has two eigenspaces $V$ and $W$ which are orthogonal and satisfy $V+W = \mathbb{R}^n$.  Hence $W=V^\perp$ (with respect to the standard inner product on $\mathbb{R}^n$) and vice versa.  Thus the existence of a nonnegative eigenvector of $M$ is an immediate corollary of the following proposition.
\begin{proposition}
For any subspace $V\subseteq \mathbb{R}^n$, either $V$ contains a nonzero, nonnegative vector or $V^\perp$ does. 
\end{proposition}

Some commentary before commencing with the proof: Although this is ostensibly a result about linear algebra, we have noted already that the notion of nonnegativity is inherently \textit{not} a purely linear algebraic property.  Hence it should not be surprising that the proof should require other ideas.  It turns out that \textit{convexity} is the key here. 

\begin{proof}[Proof of Proposition]

Define sets $$\mathbb{R}^n_{\geq 0}:=\{x\in\mathbb{R}^n : \: x_i\geq 0 \;\; \text{for all $i$}\}$$ and $$S:=\{x\in\mathbb{R}^n_{\geq 0} : \sum_i x_i = 1\}.$$  $\mathbb{R}^n_{\geq 0}$ is the set of all nonnegative vectors, and proving the proposition amounts to showing that $V$ or $V^\perp$ intersects $\mathbb{R}^n_{\geq 0}$ in a nonzero vector.  Because $V$ and $\Rp$ are both closed under multiplication by a positive scalar, if $V$ intersects $\Rp$ in a nonzero vector, then by scaling that vector appropriately we see that $V$ also intersects $S$.  So the proposition is equivalent to the statement that $V$ or $V^\perp$ intersects $S$.  

Note that $S$ is convex and compact.  A basic fact about convex sets \cite[p. 122]{separation} is the \textit{Hyperplane Separation Theorem}, which states (in one of its forms) that if two convex sets $A,B\subset \mathbb{R}^n$ are disjoint and at least one of them is compact, then there exists a constant $c\in\mathbb{R}$ and a vector $v\in \mathbb{R}^n$ such that the hyperplane defined by $x\cdot v = c$ separates $A$ and $B$.  That is, $x\cdot v > c$ for any $x\in A$, and $x\cdot v < c$ for any $x\in B$.  

If $S$ and $V$ do not intersect, then since both are convex and $S$ is compact, there exist $v$ and $c$ as above so that $x\cdot v > c$ for all $x\in S$ and $x\cdot v < c$ for all $x\in V$.  In fact, since $V$ is closed under scalar multiplication, we must have $x\cdot v = 0$ for all $x \in V$, for otherwise we could choose some scalar $\lambda \in \mathbb{R}$ so that the vector $\lambda x\in V$ satisfied $\lambda x \cdot v > c$, a contradiction.  This shows that $v\in V^\perp$ and that $c>0$.  Since all coordinate vectors $e_i$ are in $S$, all components $v_i = e_i \cdot v$ are greater than $c$ and thus positive, so that $v$ is a positive vector. 

\end{proof}

Note that the proof establishes something slightly stronger than the proposition, namely that if $V$ does not contain a nonzero, nonnegative vector, then $V^{\perp}$ must contain a \textit{positive} vector.  In general, though, neither space must contain a strictly positive vector, provided both contain nonzero, nonnegative vectors.  This is the case, for example, if $V$ is the $x$-axis and $V^\perp$ the $y$-axis in $\mathbb{R}^2$.  

The restriction to only two eigenvalues might seem to be a substantial limitation, but it is necessary.  With three or more eigenvalues, we can arrange the eigenspaces to avoid $\Rp$.  The idea is to build a matrix with two orthogonal eigenvectors $v$ and $w$ in $\mathbb{R}^n$ such that neither is itself nonnegative, but which yield a positive vector as the linear combination $v+w$.  For example, $v=(-\frac{1}{2},1,1,\dots,1)^T$ and $w=(1,-\frac{1}{2},1,0,0,0,\dots,0)^T$ are two such vectors.  (Note that this construction requires the ambient space $\mathbb{R}^n$ to have dimension $n\geq 3$, a fact that was already implicit in the hypothesis that our matrix of interest has three or more eigenvalues.)

A real symmetric matrix having $v$ and $w$ as eigenvectors with nondegenerate eigenvalues will have no nonnegative eigenvectors.  For $v$ and $w$ are not nonnegative by construction, while any other eigenvector $u$ must be orthogonal to $v$ and $w$ and thus also to $v+w$.  It is impossible for a nonzero, nonnegative vector to be orthogonal to a positive vector (the dot product would be positive), and thus $u$ cannot be nonnegative.

\begin{acknowledgment}{Acknowledgments.}
I thank the referees for helpful suggestions, including a better version of the Hyperplane Separation Theorem.  
\end{acknowledgment}

\begin{affil}
Department of Physics, Stanford University, Stanford, CA 94305 \\
hos3@cornell.edu
\end{affil}

\end{document}